\documentclass[reqno,a4paper,12pt]{amsart}
\usepackage{amssymb}
\usepackage[english]{babel}
\usepackage[dvips, lmargin=3cm, rmargin=3cm, tmargin=3cm, bmargin=3cm]{geometry}
\usepackage[colorlinks=true,urlcolor=blue,linkcolor=blue]{hyperref}
\usepackage{hyperref}
\usepackage[T1]{fontenc}
\usepackage[english]{babel}

\usepackage{dsfont}
\everymath{\displaystyle}

\usepackage{enumitem}
\usepackage{dsfont}
\newtheorem{thm}{Theorem}[section]
\newtheorem{defini}{Definition}[section]
\newtheorem{rem}{Remark}[section]
\newtheorem{lem}{Lemma}[section]
\newtheorem{prop}{Proposition}[section]
\newtheorem{coro}{Corollary}[section]

\def \R {\mathbb{R} }

\begin{document}
\title[Existence of solutions for a nonlocal Kirchhoff type problem... ]{Existence of solutions for a nonlocal Kirchhoff type problem in Fractional Orlicz-Sobolev spaces }
\author[E. Azroul, A. Benkirane, M. Srati and M.Shimi]
{E. Azroul$^1$, A. Benkirane$^2$,  M. Srati$^3$ and M. Shimi$^4$}
\address{E. Azroul, A. Benkirane,  M. Srati and M.Shimi\newline
 Sidi Mohamed Ben Abdellah
 University,
 Faculty of Sciences Dhar El Mahraz, Laboratory of Mathematical Analysis and Applications, Fez, Morocco.}
\email{$^1$elhoussine.azroul@gmail.com}
\email{$^2$abd.benkirane@gmail.com}
\email{$^3$srati93@gmail.com}
\email{$^4$mohammed.shimi2@usmba.ac.ma}

\subjclass[2010]{46E30, 58E05, 35R11, 35J60.}
\keywords{Fractional Orlicz-Sobolev spaces, Fractional $A$-Laplacian operator, Kirchhoff type problem, Direct method in variational methods.}
\maketitle
\textbf{Abstract}
In this paper, we investigate the existence of weak solution for a Kirchhoff type problem driven by a nonlocal  operator of elliptic type in a fractional Orlicz-Sobolev space, with homogeneous Dirichlet boundary conditions
 {\small$$
 (D_{K,A}) \hspace*{0.5cm} \left\{ 
   \begin{array}{clclc}

  M\left( \displaystyle \int_{\R^{2N}}A\left( [u(x)-u(y)] K(x,y)\right) dxdy\right) \mathcal{L}^K_A u & = & f(x,u)   & \text{ in }& \Omega, \\\\
   \hspace*{7cm} u & = & 0 \hspace*{0.2cm} \hspace*{0.2cm} & \text{ in } & \R^N\setminus \Omega.
    \label{eq1}
   \end{array}
   \right. 
$$ }
Where  $\mathcal{L}^K_A$  is a nonlocal operator with singular kernel $K$ and $A$ is an $N$-function, $\Omega$ is an open bounded subset in $\R^N$ with Lipschitz boundary $\partial \Omega$.
\section{Introduction}
 The  equation, 
\begin{equation}\label{1.}
\rho \dfrac{\partial^2u}{\partial t^2}-\left( \dfrac{P_0}{h}+\dfrac{E}{2L}\int_{0}^{L}|\dfrac{\partial u}{\partial x}|^2dx\right) \dfrac{\partial^2u}{\partial x^2}
\end{equation}
presented by Kirchhoff \cite{09} in 1883 is an extension of the classical d'Alembert's wave equation by considering the changes in the length of the string during vibrations. In $(\ref{1.})$, $L$ is the length of string, $h$ is the area of the cross section, $E$ is the Young modulus of the material, $\rho$ is the mass density, and $P_0$ is the initial tension. The Kirchhoff's model takes into account the length changes of the string produced by transverse vibrations. Some interesting results can be found, for example in \cite{05}. On the other hand, Kirchhoff-type boundary value problems model several physical and biological systems where $u$ describes a process which depend on the average of itself, as for example, the population density. We refer the reader to \cite{03,04,013} for some related works. In \cite{06}, the authors showed the existence and multiplicity of solutions to a class of $p(x$)-Kirchhoff type equations via variational methods. In \cite{09}, the author obtained the existence of infinite solutions to the $p$-Kirchhoff type quasilinear elliptic equations via the fountain theorem. In \cite{006}, the authors investigated
higher order $p(x)$-Kirchhoff type equations via symmetric mountain pass Theorem, even in the degenerate case. However, they did not consider the existence of solutions for Kirchhoff type problems in the fractional setting. In the very recent paper \cite{08}, the authors first provided a detailed discussion about the physical meaning underlying the fractional Kirchhoff models and their applications, see \cite{08} for further details.

 Problems of this type have been intensively studied in the last few years, due to numerous and relevant applications in many fields of mathematics, such as approximation theory, mathematical physics (electrorheological fluids), calculus of variations, nonlinear potential theory, the theory of quasiconformalmappings, differential geometry, geometric function theory, probability theory and image processing (see for instance \cite{012,06.,04.}).

 The study of nonlinear elliptic equations involving quasilinear homogeneous type operators is based on the theory of Sobolev spaces and fractional Sobolov spaces $W^{s,p}(\Omega)$  in order to find weak solutions. In certain equations, precisely in the case of nonhomogeneous differential operators, when trying to relax some conditions on these operators (as growth conditions), the problem can not be formulated with classical Lebesgue and Sobolev spaces. Hence, the adequate functional spaces is the so-called Orlicz spaces. These spaces consists of functions that have weak derivatives and satisfy certain integrability conditions. Many properties of Orlicz-Sobolev spaces and fractional Orlicz-Sobolev spaces come in \cite{1,3,7,12,13}. For this, many researchers have studied the existence of solutions for the eigenvalue problems involving nonhomogeneous operators in the divergence form through Orlicz-Sobolev spaces by using variational methods and critical point theory, monotone operator methods, fixed point theory and degree theory (see for instance \cite{01,001,0001,07.,007.,0007.,00007.,012.,5,10.,37.}).

We consider the following problem, 
{\small$$
 (D_{K,p}) \hspace*{0.5cm} \left\{ 
   \begin{array}{clclc}

  M\left( \displaystyle \int_{\R^{2N}}|u(x)-u(y)|^p K(x-y) dxdy\right) \mathcal{L}^p_K u & = & f(x,u)   & \text{ in }& \Omega \\\\
   \hspace*{7cm} u & = & 0 \hspace*{0.2cm} \hspace*{0.2cm} & \text{ in } & \R^N\setminus \Omega,
    \label{eq1}
   \end{array}
   \right. 
$$ }
where $N>sp$ with $s\in (0,1)$,  $\Omega$ is an open bounded subset in $\R^N$ with Lipschitz boundary $\partial \Omega$, $M: [0,\infty) \longrightarrow (0,\infty)$ is a continuous function, $f: \Omega\times \R \longrightarrow \R$ is a Carath\'eodory function and $\mathcal{L}^p_K u$ is a non-local operator
defined as follows:
{\small $$
        \begin{aligned}
        \mathcal{L}^p_K u(x)&=2P.V \int_{\R^N} | u(x)-u(y)|^{p-2}(u(x)-u(y))K(x-y)dy,\hspace{0.2cm} x\in \R^N,\\
        &=2\lim\limits_{\varepsilon\searrow 0}\int_{\R^N \smallsetminus B_\varepsilon(x)}  | u(x)-u(y)|^{p-2}(u(x)-u(y))K(x-y)dy,\hspace{0.2cm} x\in \R^N,
        \end{aligned}
         $$}
where $1 < p < \infty$ and $K : \R^N \setminus \left\lbrace 0\right\rbrace \longrightarrow (0,\infty)$ is a measurable function with the following property :
{\small$$
 \left\{ 
   \begin{array}{clclc}

  & K(x)= K(-x) &   & \text{ for all }& x\in \R^N\setminus \left\lbrace 0\right\rbrace, \\\\
    &\text{ there exist } k_0>0  \text{ such that }&  K(x)\geqslant \dfrac{k_0}{|x|^{N+sp}} &\text{ for all }& x\in \R^N\setminus \left\lbrace 0\right\rbrace,\\\\
     &\delta K\in L^1(\R^N)&  &\text{ where }& \delta(x)=\min\left\lbrace 1,|x|^p \right\rbrace .
    \label{eq1}
   \end{array}
   \right. 
$$ }

 In \cite{111}, the authors discuss the above-mentioned problem in two cases: when $f$ satisfies sublinear
 growth condition, the existence of nontrivial weak solutions is obtained by applying
 the direct method in variational methods; when $f$ satisfies suplinear growth
 condition, the existence of two nontrivial weak solutions is obtained by using the
 mountain pass Theorem.

In this paper, our purpose is to generalize the previous works in the setting of fractional Orlicz-Sobolev spaces. That is we are interested to study the existence of weak solution for a Kirchhof type problem driven by a nonlocal operator of elliptic type, with homogeneous Dirichlet boundary conditions as follows:

 {\small$$
 (D_{K,A}) \hspace*{0.5cm} \left\{ 
   \begin{array}{clclc}

  M\left( \displaystyle \int_{\R^{2N}}A\left( [u(x)-u(y)] K(x,y)\right) dxdy\right) \mathcal{L}^K_A u & = & f(x,u)   & \text{ in }& \Omega \\\\
   \hspace*{7cm} u & = & 0 \hspace*{0.2cm} \hspace*{0.2cm} & \text{ in } & \R^N\setminus \Omega,
    \label{eq1}
   \end{array}
   \right. 
$$ }

where $N>1$, $\Omega$ is an open bounded subset in $\R^N$ with Lipschitz boundary $\partial \Omega$, $M: [0,\infty) \longrightarrow (0,\infty)$ is a continuous function $f: \Omega\times \R \longrightarrow \R$ is a Carath\'eodory function satisfying the Ambrosetti-Rabinowitz type condition and $\mathcal{L}^K_A u$ is a nonlocal integro-differential operator of elliptic type defined as follows:

 {\small $$
        \begin{aligned}
        \mathcal{L}^K_A u(x)&=2P.V \int_{\R^N} A'\left( (u(x)-u(y))K(x,y)\right)K(x,y)dy, \hspace{0.2cm} x\in \R^N,\\
        &=2\lim\limits_{\varepsilon\searrow 0} \int_{\R^N\setminus B_\varepsilon(x)} a\left((u(x)-u(y))K(x,y)\right)K(x,y)dy,\hspace{0.2cm} x\in \R^N,
        \end{aligned}
         $$}
 with $A$ is an $N$-function and the singular kernel $K: \R^N\times \R^N \longrightarrow (0,\infty)$ is a measurable function satisfies the following  properties :

{\small$$
 (K_0) \left\{ 
   \begin{array}{clclc}

   & K(x,y)=K(y,x)&    & \text{ for all }& (x,y)\in \R^N\times\R^N,  \\\\
    & \exists k_0>0  \text{ : }   K(x,y)\geqslant \dfrac{k_0}{|x-y|^sA^{-1}(|x-y|^N)}& &\text{ for all }& (x,y)\in \R^N\times\R^N,\\\\
    & \delta K\in L_M(\R^N\times\R^N) & &\text{ where }& \delta(x,y)=\min\left\lbrace 1,|x-y| \right\rbrace.
    \label{eq1}
   \end{array}
   \right. 
$$ }
A typical example for $K$ is given by singular kernel $K(x,y)=\frac{1}{|x-y|^sA^{-1}(|x-y|^N)}$. In this case $\mathcal{L}^K_A=(-\Delta)^s_A$ which is introduced in \cite{3} and the problem $(D_{K,A})$ becomes

 {\small$$
(D_{s,A})  \left\{ 
   \begin{array}{clclc}

  M\left( \displaystyle \int_{\R^{2N}}A\left(\dfrac{ u(x)- u(y)}{|x-y|^sA^{-1}(|x-y|^N)}\right) dxdy\right) (-\Delta)^s_A u & = & f(x,u)   & \text{ in }& \Omega \\\\
   \hspace*{7cm} u & = & 0 \hspace*{0.2cm} \hspace*{0.2cm} & \text{ in } & \R^N\smallsetminus \Omega,
    \label{eq1}
   \end{array}
   \right. 
$$ }
where $(-\Delta)^s_A$ is the fractional $A$-Laplacian operator  defined on smooth functions by :
 {\small $$
        \begin{aligned}
        (-\Delta)^s_Au(x)&=2P.V \int_{\R^N} A'\left( \dfrac{(u(x)-u(y))}{|x-y|^sA^{-1}(|x-y|^{N}) }\right)\dfrac{dy}{|x-y|^s A^{-1}(|x-y|^{N})}\\
        &=2\lim\limits_{\varepsilon\searrow 0} \int_{\R^N\setminus B_\varepsilon(x)} a\left( \dfrac{(u(x)-u(y))}{|x-y|^sA^{-1}(|x-y|^{N}) }\right) \dfrac{dy}{|x-y|^sA^{-1}(|x-y|^{N}) }
        \end{aligned}
         $$}
Note that for $A(t)=\frac{|t|^p}{p}$, the problem $(D_{K,A})$  co\"incide with the problem $(D_{K,p})$.    

This paper is organized as follows : in the second section, we recall some properties of Orlicz-Sobolev and fractional Orlicz-Sobolev spaces into Lebesgue  spaces. In the third section, we will introduce some necessary definitions and properties of space $W^KL_A$. Moreover, we prove the continuous and compact embedding theorems of these spaces. 
Finally,   using the direct method in calculus variations, we obtain the existence of a weak solution of problem  $(D_{K,A})$.

\section{Some preliminaries results}

To deal with this situation we introduce the fractional Orlicz-Sobolev space to investigate problem $(D_{K,A})$. Let us recall the definitions and some elementary properties of these spaces. For more details  we refer the reader to \cite{1,3,6,11,12}.\\

 Let $\Omega$ be an open subset of $\R^N$. Let $A$ :
  $\R^+ \rightarrow \R^+$ be an N-function, that is, $A$ is continuous, convex, with $A(t) > 0$ for $t >
  0$, $ \frac{A(t)}{t}\rightarrow 0$ as  $t \rightarrow 0$ and $\frac{A(t)}{t} \rightarrow\infty$ as $t \rightarrow \infty$. Equivalently, $A$ admits the
  representation : $ A(t)=\displaystyle\int_{0}^{t} a(s)ds$ where $a : \R^+ \rightarrow \R^+$ is non-decreasing, right
  continuous, with $a(0) = 0$, $ a(t) > 0$, for all $t> 0$ and $a(t) \rightarrow \infty $ as $t \rightarrow \infty$.  The conjugate N-function of $A$ is defined by $\overline{A}(t) =\displaystyle \int_{0}^{t} \overline{a}(s)ds$, where $ \overline{a} : \R^+\rightarrow
  \R^+$ is given by $\overline{a}(t) = \sup \left\lbrace s : a(s) \leqslant t\right\rbrace$.\\ Evidently we have
  \begin{equation}\label{1}
st\leqslant A(t)+\overline{A}(s),
  \end{equation}
  which is known Young's inequality. Equality holds in (\ref{1}) if and only if either $t=\overline{a}(s)$ or $s=a(t)$.
 
  We say that the N-function $A$  satisfy the
  global $\Delta_2$-condition  if, for some $k > 0,$
 $$ A(2t) \leqslant kA(t)\text{ , }\forall t \geqslant 0.$$
 When this inequality holds only for $t \geqslant t_0 > 0$, $M$ is said to satisfy the $\Delta_2$-condition
   near infinity.
 \\
 We call the pair $(A,\Omega)$ is $\Delta$-regular if either :\\
 (a) $A$  satisfies a
   global $\Delta_2$-condition, or \\
   (b) $A$ satisfies a $\Delta_2$-condition near infinity and $\Omega$ has finite volume.
 
 According to \cite{1}, $A$ satisfies the
   global $\Delta_2$-condition  if and only if 
  $$ (a_1)   \hspace{4cm} 1\leqslant \inf_{s\geqslant 0}\dfrac{sa(s)}{A(s)}<\sup_{s\geqslant 0}\dfrac{sa(s)}{A(s)}<+\infty. \hspace{5cm} \text{.}$$
 We assume that : 
 $$ (a_2)   \hspace{3cm} \text{ the function } t\mapsto A(\sqrt{t}) \text{ is convex for all } t\geqslant 0 . \hspace{5cm} \text{.}$$
 
    The Orlicz class $K_A (\Omega)$ (resp. the Orlicz space
   $L_A(\Omega)$) is defined as the set of (equivalence classes of) real-valued measurable
   functions $u$ on $\Omega$ such that  
\begin{equation}\label{3}
\int_{\Omega} A(u(x))dx <\infty \hspace*{0.5cm} \text{ (resp. } \int_{\Omega} A(\lambda u(x))dx< \infty \text{ for some } \lambda >0 ).
\end{equation}
$L_A(\Omega)$ is a Banach space under the Luxemburg-norm
\begin{equation}\label{4}
 ||u||_A=\inf \Bigg\{\lambda>0  : \int_{\Omega}A\left(  \dfrac{u(x)}{\lambda}\right) dx \leqslant 1\Bigg\},
\end{equation}
and $K_A(\Omega)$ is a convex subset of $L_A(\Omega)$. The closure in $L_A(\Omega)$ of the set of
bounded measurable functions on $\Omega$ with compact support in $\overline{\Omega}$ is denoted by $E_A(\Omega)$.\\
The equality $E_A(\Omega)=L_A(\Omega)$ holds if and only if $A$ is satisfies the global $\Delta_2$-condition.

 Using the Young's inequality, it is possible to prove a H\"older type inequality,
  that is,
  \begin{equation}
  \left| \int_{\Omega}uvdx\right| \leqslant 2||u||_A||v||_{\overline{A}}\hspace*{0.5cm} \text{  } \forall u \in L_A(\Omega)  \text{ and } v\in L_{\overline{A}}(\Omega).
  \end{equation}
  
Now, we defined the fractional Orlicz-Sobolev space $W^sL_A(\Omega)$ as follows 
\begingroup\makeatletter\def\f@size{10}\check@mathfonts
$$W^s{L_A}(\Omega)=\Bigg\{u\in L_A(\Omega) : \int_{\Omega} \int_{\Omega} A\left( \dfrac{\lambda( u(x)- u(y))}{|x-y|^sA^{-1}(|x-y|^N)}\right) dxdy< \infty, \text{ for some } \lambda>0 \Bigg\}.
$$\endgroup
This space is equipped with the norm,
\begin{equation}\label{6}
||u||_{s,A}=||u||_{A}+[u]_{s,A},
\end{equation}

where $[.]_{s,A}$ is the Gagliardo seminorm, defined by 
$$[u]_{s,A}=\inf \Bigg\{\lambda >0 :  \int_{\Omega} \int_{\Omega} A\left( \dfrac{u(x)- u(y)}{\lambda|x-y|^sA^{-1}(|x-y|^N)}\right) dxdy\leqslant 1 \Bigg\}.
$$
We define $W^s_0L_A(\Omega)$ as the set of the functions $u$ in  $W^sL_A(\Omega)$ such that $u=0$ in $\R^N\setminus \Omega$. 
 By \cite{3}, $W^sL_A(\Omega)$ and $W^s_0L_A(\Omega)$ are separable Banach spaces, (resp. reflexive spaces)  if and only if $A$ satisfies the global $\Delta_2$-condition (resp. $A$ and  $\overline{A}$ are satisfy the global $\Delta_2$-condition).

Let $A$ be a given N-function, satisfying the following conditions :
  \begin{equation}\label{5}
  \int_{0}^{1} \dfrac{A^{-1}(\tau)}{\tau^{\frac{N+s}{N}}}d\tau<\infty,
  \end{equation}
  
  \begin{equation}\label{6.}
  \int_{1}^{\infty} \dfrac{A^{-1}(\tau)}{\tau^{\frac{N+s}{N}}}d\tau=\infty.
  \end{equation}
  For instance if $A(t)=\frac{1}{p}t^p$, then (\ref{5}) holds precisely when $sp<N$.\\
  If (\ref{6.}) is satisfied, we define the inverse Sobolev conjugate N-function of $A$ as follows, 
  \begin{equation}\label{7}
  A_*^{-1}(t)=\int_{0}^{t}\dfrac{A^{-1}(\tau)}{\tau^{\frac{N+s}{N}}}d\tau.
  \end{equation}
\begin{thm}\label{2.1} (cf. \cite{3})
 Let $A$ be an $N$-function and $s\in (0,1)$. Let $\Omega$  be a bounded open
  subset of  $\R^N$ with $C^{0,1}$-regularity 
    and bounded boundary. If $(\ref{5})$ and $(\ref{6.})$  hold, then 
 \begin{equation}\label{18}
  W^s{L_A}(\Omega)\hookrightarrow L_{A_*}(\Omega).
 \end{equation}
\end{thm}

 \begin{thm}\label{2.2.}(cf. \cite{3})
        Let $s\in (0,1)$ and $A$ be an $N$-function. Let $\Omega$  be a bounded open
          subset of  $\R^N$ and  $C^{0,1}$-regularity 
            with bounded boundary. If $(\ref{5})$ and  $(\ref{6.})$ hold, then the embedding
         \begin{equation}\label{27}
          W^s{L_A}(\Omega)\hookrightarrow L_{B}(\Omega),
         \end{equation}
         is compact for all $B\prec\prec A_*$.
         \end{thm}  
         
  By fixing the fractional exponent $s\in(0,1)$ and for any $p\in [1,\infty)$, we define the fractional Sobolev space $W^{s,p}(\Omega)$ as follows,
  $$ W^{s,p}(\Omega)=\Bigg\{u \in L^p(\Omega)  \text{ : }  \dfrac{|u(x)-u(y)|}{|x-y|^{\frac{N}{p}+s}} \in L^p(\Omega \times \Omega)\Bigg\};$$
  that is, an intermediary Banach space between, endowed
  with the natural norm
  $$||u||_{s,p}=\Bigg(\int_{\Omega}|u|^pdx+\int_{\Omega}\int_{\Omega} \dfrac{|u(x)-u(y)|^p}{|x-y|^{sp+N}}dxdy\Bigg)^{\frac{1}{p}}.$$
  \begin{thm}$\label{2.3}$(cf. \cite{11}).
   Let $s\in (0,1)$ and let $p\in [1,+\infty)$ such that $sp<N$. Let $\Omega$  be an open
    subset of  $\R^N$  with $C^{0,1}$-regularity and bounded boundary. So there exists a positive constant $C=C(N,s,p,\Omega)$ such that, for all  $f \in W^{s,p}(\Omega)$ we have  \\
     $$||f||_{L^q(\Omega)} \leqslant C||f||_{W^{s,p}(\Omega)} \text{    } \forall q \in [p,p^*],$$
     that is,
     $$ W^{s,p}(\Omega)  \hookrightarrow L^q(\Omega)  \text{    } \forall q \in[p,p^*],$$
       where $p^*=\frac{Np}{N-sp}$.\\
       If, in addition, $\Omega$ is bounded, then the space $W^{s,p}(\Omega)$ is continuously
       embedded in $L^q(\Omega)$ for any $q\in [1,p^*]$.
  \end{thm}
  \begin{thm}$\label{2.4}$(cf. \cite{6}).
   Let $s\in (0,1)$ and let $p\in [1,+\infty)$ such that $sp<N$. Let $\Omega$  be a bounded open
    subset of  $\R^N$   with $C^{0,1}$-regularity and bounded boundary. Then the embedding 
   $$ W^{s,p}(\Omega)  \hookrightarrow L^q(\Omega)  \text{    } \forall q \in[1,p^*),$$
   is compact.
  \end{thm}

      \begin{thm}\label{2.2} (cf. \cite{110})
     Suppose $X$ is a reflexive Banach space with norm $||.||$ and let
     $V\subset X$ be a weakly closed subset of $X$. Suppose $E : V \longrightarrow \R \cup \left\lbrace +\infty\right\rbrace $ is coercive
     and (sequentially) weakly lower semi-continuous on $V$ with respect to $X$ , that
     is, suppose the following conditions are fulfilled:
     
     (1) $E(u)\rightarrow \infty$ as $||u||\rightarrow \infty$, $u\in V$.
     
     (2)  For any $u\in V$, any sequence $\left\lbrace u_n\right\rbrace $ in $V$ such that $u_n\rightharpoonup u$ weakly in $X$
     there holds:
    $$E(u)\leqslant \liminf_{n\rightarrow \infty}E(u_n).$$
     Then $E$ is bounded from below on $V$ and attains its infimum in $V$.
     
      \end{thm}
  \section{Variational framework}
  One of the aims of this paper is to study nonlocal problems driven by $(-\Delta)^s_A$
  (or its generalization $\mathcal{L}^K_A$)  with Dirichlet boundary data via variational methods. For this purpose, we need to work in a suitable fractional Orlicz-Sobolev space: for
  this, we consider a functional analytical setting that is inspired by (but not equivalent to) the fractional Orlicz-Sobolev spaces in order to correctly encode the Dirichlet
  boundary datum in the variational formulation.
  
In this section, we show some basic results that will be used in the next section.
\begin{defini}
Let $\Omega$ be an open subset of $\R^N$, we denote $Q=\R^N\setminus \mathcal{O}$, where $\mathcal{O}=(\R^N\setminus \Omega)\times (\R^N\setminus\Omega)$, and let $A$ be an $N$-function and $K$ as defined in $(K_0)$. We define the space $W^KL_A(\Omega)$ since the set of Lebesgue measurable functions from $\R^N$ to $\R$ such that $ u_{|_{\Omega}}$ in $L_A(\Omega)$ and 
$$ (u(x)-u(y))K(x,y) \in L_A(Q).$$

\end{defini}
The space $W^KL_A$ is equipped with the norm
\begin{equation}
||u||_{K,A}=||u||_A+[u]_{K,A}
\end{equation}
where $[.]_{K,A}$ is the Gagliardo seminorm, defined by :
$$[u]_{K,A}=\inf \Bigg\{\lambda >0 :  \int_{Q}  A\left( \dfrac{u(x)- u(y)}{\lambda}K(x,y)\right) dxdy\leqslant 1 \Bigg\}.
$$
$||u||_{K,A}$ is a norm in $W^KL_A(\Omega)$. Indeed, let $u\in W^KL_A(\Omega)$ such that  $u=0$ then for all $\lambda >0$, we get
$$||u||_A=0 \text{ and } \int_Q A\left( \dfrac{u(x)- u(y)}{\lambda}K(x,y)\right)dxdy=0,$$
this implies that 
$$ ||u||_A=0 \text{ and } [u]_{K,A}=0. $$
Conversely, if $||u||_{K,A}=0$, this implies that 
$$u=0 \text{ a.e in } \Omega \text{ and } \int_{Q}  A\left( \dfrac{u(x)- u(y)}{\lambda}K(x,y)\right) dxdy\leqslant 1$$
for all $\lambda>0$. Let $0<h\leqslant 1$ then
$$\int_{Q}  A\left(\frac{h}{h} (u(x)- u(y))K(x,y)\right) dxdy  \leqslant h \int_{Q}  A\left( \frac{1}{h}(u(x)- u(y))K(x,y)\right) dxdy\leqslant h$$
so for $h\rightarrow 0^+$ we have
$$\int_{Q}  A\left( (u(x)- u(y))K(x,y)\right) dxdy=0$$ 
then $u(x)=u(y)$ a.e in $Q$ therefore $u=c\in \R$ a.e. in $\R^N$, since $u=0$ a.e in $\Omega$, we get $u=0$ in $\R^N$.

We shall work in the closed linear subspace

$$W^K_0L_A(\Omega)=\left\lbrace u\in W^KL_A(\Omega) \text{ : } u=0 \text{ in } \R^N\setminus \Omega\right\rbrace. $$

\begin{lem} \label{3.1..}
$C^2_0(\Omega)\subset W^K_0L_A(\Omega)$.
\end{lem}
      \noindent \textbf{Proof}. 
Let $u\in C^2_0(\Omega)$, since $u$ vanishes outside $\Omega$, we only need that to check that 
$$ \int_{\R^N}\int_{\R^N}  A\left(\lambda (u(x)- u(y))K(x,y)\right) dxdy <\infty \text{ for some } \lambda>0.$$
Indeed, 
$$
 \begin{aligned}
\int_{\R^N}\int_{\R^N}  A\left( (u(x)- u(y))K(x,y)\right) dxdy &= \int_{\Omega}\int_{\Omega}  A\left( (u(x)- u(y))K(x,y)\right) dxdy\\ &+2\int_{\Omega}\int_{\R^N \setminus \Omega}  A\left((u(x)- u(y))K(x,y)\right) dxdy\\
&\leqslant 2\int_{\Omega}\int_{\R^N}  A\left((u(x)- u(y))K(x,y)\right) dxdy.
\end{aligned} 
$$
Now we notice that
$$|u(x)-u(y)|\leqslant ||\nabla u||_{L^{\infty}(\R^N)}|x-y| \text{ and } |u(x)-u(y)|\leqslant 2||u||_{L^{\infty}(\R^N)}.$$
Accordingly, we get 
$$|u(x)-u(y)|\leqslant 2||u||_{C^{1}(\R^N)}\min\left\lbrace 1,|x-y|\right\rbrace:=\alpha \delta(x,y). $$
with $\alpha =2||u||_{C^{1}(\R^N)}$ and since $\delta K \in L_A(\R^N \times \R^N)$. There exist $\lambda>0$, such that,
$$\int_{\R^N}\int_{\R^N}  A\left(\dfrac{\lambda}{ \alpha}(u(x)- u(y))K(x,y)\right) dxdy\leqslant \int_{\R^N}\int_{\R^N}  A\left( \lambda \delta(x,y)K(x,y)\right) dxdy<\infty$$
this implies that $(u(x)-u(y))K(x,y)\in L_A(\R^N \times \R^N)$.\\
                    \hspace*{15cm$\Box$ }  
\begin{rem}
A trivial consequence of lemma $\ref{3.1..}$, $W^KL_A(\Omega)$ and $W^K_AL(\Omega)$ are non-empty.
\end{rem}                    
  \begin{lem}\label{3.2}
  Let $s\in (0,1)$ and let $K : \R^N \times \R^N \rightarrow (0,\infty)$ satisfy assumption $(K_0)$, then the following assertion hold :
  
  1) $W^KL_A(\Omega) \hookrightarrow W^sL_A(\Omega)$,
  
    2) $W_0^KL_A(\Omega) \hookrightarrow W^sL_A(\R^N)$.
  \end{lem}
        \noindent \textbf{Proof}.
        
      1)  Let $u\in W^KL_A(\Omega)$, and $\lambda>0$, by $(K_0)$ we get 
       $$ \int_{\Omega}\int_{\Omega}A\left(\frac{k_0}{\lambda} \dfrac{u(x)-u(y)}{|x-y|^sA^{-1}(|x-y|^{N})}\right) dxdy\leqslant \int_{Q}  A\left( \dfrac{u(x)- u(y)}{\lambda}K(x,y)\right) dxdy,$$
         then $u\in W^sL_A(\Omega)$ and 
     $$[u]_{s,A}\leqslant \frac{1}{k_0}[u]_{K,A},$$
     this implies that 
     $$||u||_{s,A}\leqslant \sup\left\lbrace 1,\frac{1}{k_0}\right\rbrace ||u||_{K,A}.
  $$
  
  2) Let $u \in W^K_0L_A(\Omega)$, so $u=0$ in $\R^N\setminus \Omega$ and 
  $$||u||_{L_A(\R^N)}=||u||_{L_A(\Omega)}<\infty.$$
  On the other hand, for all $\lambda>0$ we get
  \begingroup\makeatletter\def\f@size{10.5}\check@mathfonts
   $$
   \begin{aligned}
    \int_{\R^N}\int_{\R^N}A\left(\frac{k_0}{\lambda} \dfrac{u(x)-u(y)}{|x-y|^sA^{-1}(|x-y|^{N})}\right) dxdy&=\int_{Q}A\left(\frac{k_0}{\lambda} \dfrac{u(x)-u(y)}{|x-y|^sA^{-1}(|x-y|^{N})}\right) dxdy\\
    &\leqslant   \int_{Q}  A\left( \dfrac{u(x)- u(y)}{\lambda}K(x,y)\right) dxdy,
    \end{aligned}
    $$\endgroup
    then $u\in W^sL_A(\R^N)$ and 
         $$[u]_{W^sL_A(\R^N)}\leqslant \frac{1}{k_0}[u]_{K,A}$$
      this implies that 
           $$||u||_{W^sL_A(\R^N)}\leqslant \sup\left\lbrace 1,\frac{1}{k_0}\right\rbrace ||u||_{K,A}.
        $$
                                       \hspace*{15cm$\Box$ }  
\begin{thm}\label{3.1}(Generalized Poincar\'e inequality). 
Let $\Omega$ be a bounded open subset of $\R^N$, let $K$ as defined in $(K_0)$, and let $A$ be an $N$-function. Then there exist a positive constant $\mu$ such that 
$$||u||_{A}\leqslant \mu [u]_{K,A} \text{ for all } u\in W^K_0L_M(\Omega).$$
\end{thm} 
    To prove Theorem $\ref{3.1}$, we need to the following lemma.
   
       \begin{lem}\label{3.3}
       Let $\Omega$ be a bounded open subset of  $\R^N,$ and let $s\in (0,1)$. Let $A$ be an N-function. Then there exists a positive constants $\alpha$ and $\beta$ such that, \\  
                \begin{equation}\label{11}
                  \int_{\Omega}A\left( \dfrac{u}{\alpha \lambda}\right) dx \leqslant \beta \int_{\Omega}\int_{\Omega}A\left( \dfrac{u(x)-u(y)}{\lambda|x-y|^sA^{-1}(|x-y|^{N})}\right) dxdy,
                  \end{equation}
           for all $u \in W^s_0L_A(\Omega)$ and all $\lambda>0$. In particular 
     $$||u||_{A}\leqslant \mu [u]_{s,A} \text{ for all } u\in W_0^sL_A(\Omega)$$   
     where $\mu = \alpha \beta$.          
        \end{lem}
  \noindent \textbf{Proof of Lemma \ref{3.3}}. 
                  Let $ u \in W^s_0L_A(\Omega)$ and $ B_R \subset \R ^ N \setminus \Omega $, i.e, the ball of radius $ R>0 $ in the complement of $ \Omega $. Then for all $x\in \Omega$, $y\in B_R$ and all $\lambda>0$ we have, \\
                 $$A(u(x))=A\left( \dfrac{u(x)-u(y)}{ |x-y|^s A^{-1}(|x-y|^N)}|x-y|^sA^{-1}(|x-y|^N)\right),$$
                 this implies that,\\
                 $$A(u(x))\leqslant A\left( \dfrac{u(x)-u(y)}{ |x-y|^sA^{-1}(|x-y|^N)} diam(\Omega\cup B_R)^s (A^{-1}(diam(\Omega\cup B_R)^N))\right),$$
       we suppose $\alpha =diam(\Omega\cup B_R)^s (A^{-1}(diam(\Omega \cup B_R)^N))$, we get

                {\small $$|B_R|A\left( \dfrac{u(x)}{\alpha \lambda}\right) \leqslant  \int_{B_R} A\left( \dfrac{u(x)-u(y)}{ \lambda|x-y|^sA^{-1}(|x-y|^N)} \right) dy,$$}
       then
                  $$ \int_{\Omega} A(\dfrac{u(x)}{\alpha\lambda})dx \leqslant \beta\int_{\Omega}\int_{\Omega} A\left( \dfrac{u(x)-u(y)}{ \lambda|x-y|^sA^{-1}(|x-y|^N)} \right)dxdy,$$ 
                            where $\beta =\dfrac{1}{|B_R|}$. 
                 Therefore, for $\mu=\alpha\beta$ we get
        $$||u||_{A}\leqslant \mu [u]_{s,A} \text{ for all } u\in W_0^sL_M(\Omega).$$ 
                      \hspace*{15cm$\Box$ }  
  \noindent \textbf{Proof of Theorem \ref{3.1}}.  Let $u\in W^K_0L_A(\Omega)$, by Lemma \ref{3.2} we get $u\in W^s_0L_A(\Omega)$, so by Lemma \ref{3.3}, there exists $\alpha, \beta>0$ such that 
  $$ \int_{\Omega} A\left( \dfrac{u(x)}{\alpha\lambda}\right) dx \leqslant \beta\int_{\Omega}\int_{\Omega} A\left( \dfrac{u(x)-u(y)}{ \lambda|x-y|^sA^{-1}(|x-y|^N)} \right)dxdy,$$ 
  for all $\lambda>0$. On the other hand by condition $(K_0)$, we have
  $$\int_{\Omega}A\left( \dfrac{k_0 u(x)}{\alpha\lambda}\right) dx \leqslant \beta \int_{Q}  A\left( \dfrac{u(x)- u(y)}{\lambda}K(x,y)\right) dxdy,$$
  this implies that 
  $$||u||_A\leqslant \mu[u]_{K,A}$$
  where $\mu=\dfrac{\alpha\beta}{k_0}$.\\
                      \hspace*{15cm$\Box$ } 
                      
   As a consequence of Theorem $\ref{3.1}$ $ [.] _ {K, M}$ is a norm on  $W^K_{0}L_M (\Omega) $ equivalent to $ ||. ||_{K,M}.$ Then $\left(W^K_0L_M , [.]_{K,A} \right)$ is Banach space.  
                     
  Now,   we introduce the following notations : 
               \begin{equation}\label{12}
                p_0:=\inf_{s\geqslant 0}\dfrac{sa(s)}{A(s)} \text{  , }
               \hspace{1cm}  
                 p^0:=\sup_{s\geqslant 0}\dfrac{sa(s)}{A(s)}.
                 \end{equation}                  
\begin{prop}\label{pro3}
    Assume that condition $(a_1)$ is satisfied. Then the following relations hods true 
    \begin{equation}\label{32}
    [u]^{p_0}_{K,A}\leqslant \phi(u) \leqslant [u]^{p^0}_{K,A} \text{   ,  } \forall u\in W^K_0L_A(\Omega) \text{ with }[u]_{K,A}>1,
    \end{equation} 
     \begin{equation}\label{33}
       [u]^{p^0}_{K,A}\leqslant \phi(u) \leqslant [u]^{p_0}_{K,A} \text{   ,  } \forall u\in W^K_0L_A(\Omega) \text{ with }[u]_{K,A}<1.
       \end{equation} 
    where $\phi(u):=\displaystyle\int_{Q}A\left( (u(x)-u(y))K(x,y)\right) dxdy$.
    \end{prop}
    \noindent \textbf{Proof}. First we show that  $\phi(u) \leqslant [u]^{p^0}_{K,A} \text{   for all  }  u\in W^K_0L_A(\Omega) \text{ with }[u]_{K,A}>1$. Indeed, since $p^0\geqslant \dfrac{ta(t)}{A(t)} $ for all $t>0$ it follows that for all $\sigma >1,$ we have 
    $$ \log(A(\sigma t))-\log(A(t))=\int_{t}^{\sigma t} \dfrac{a(\tau)}{A(\tau)}d\tau   \leqslant \int_{t}^{\sigma t} \dfrac{p^0}{\tau}d\tau=\log(\sigma^{p^0}).$$ 
    Thus, we deduce 
    \begin{equation}\label{34}
    A(\sigma t)\leqslant \sigma^{p^0}A(t) \text{ for all } t>0 \text{ and } \sigma>1.
    \end{equation}
    Let now $u\in W_0^KL_A(\Omega)$ with $[u]_{K,A}>1$. Using the definition of  Luxemburg norm and the relation $(\ref{34})$, we deduce 
    $$
    \begin{aligned}
    \displaystyle\int_{Q}A\left( (u(x)-u(y))K(x,y)\right) dxdy&=\displaystyle\int_{Q}A\left([u]_{K,A} \dfrac{(u(x)-u(y))K(x,y)}{[u]_{K,A}} \right) dxdy\\
    &\leqslant [u]^{p^0}_{K,A}\displaystyle\int_{Q}A\left( \dfrac{(u(x)-u(y))K(x,y)}{[u]_{K,A}} \right) dxdy\\
    &\leqslant [u]^{p^0}_{K,A}.
    \end{aligned}
    $$
    Now, we show that $\phi(u) \geqslant [u]^{p_0}_{K,A} \text{   for all  } u\in W^K_0L_A(\Omega) \text{ with }[u]_{K,A}>1$. Indeed, since $$p_0\leqslant \dfrac{ta(t)}{A(t)}$$ for all $t\geqslant 0$, it follows that for all $\sigma>1$, we have
     $$ \log(A(\sigma t))-\log(A(t))=\int_{t}^{\sigma t} \dfrac{a(\tau)}{A(\tau)}d\tau   \geqslant \int_{t}^{\sigma t} \dfrac{p_0}{\tau}d\tau=\log(\sigma^{p_0}).$$ 
      Hence, we deduce 
        \begin{equation}\label{35}
        A(\sigma t)\geqslant \sigma^{p_0}A(t) \text{ for all } t>0 \text{ and } \sigma>1.
        \end{equation}
   Let $u\in W_0^KL_A(\Omega)$ with $[u]_{K,A}>1$, we consider $\beta\in (1,[u]_{K,A})$, since $\beta<[u]_{K,A}$, so by definition of Luxemburg norm, it follows that 
 $$\int_{Q}A\left( \dfrac{(u(x)-u(y))K(x,y)}{\beta} \right) dxdy>1,$$
 the above inequality implies that 
 $$
 \begin{aligned}
 \displaystyle\int_{Q}A\left( (u(x)-u(y))K(x,y)\right) dxdy&=\int_{Q}A\left( \beta\dfrac{(u(x)-u(y))K(x,y)}{\beta} \right) dxdy\\
 &\geqslant \beta^{p_0}\int_{Q}A\left( \dfrac{(u(x)-u(y))K(x,y)}{\beta} \right) dxdy\\
 &\geqslant \beta^{p_0},
 \end{aligned}
 $$
 letting $\beta\nearrow [u]_{K,A}$, we deduce that relation $(\ref{32})$ hold true.
 
 Next, we show that $\phi(u) \leqslant [u]^{p_0}_{K,A} \text{   for all  }  u\in W^K_0L_A(\Omega) \text{ with }[u]_{K,A}<1$. By the same argument in the proof of  $(\ref{34})$ and $(\ref{35})$, we have 
 \begin{equation}\label{36}
 A(t)\leqslant \tau^{p_0}A\left( \dfrac{t}{\tau} \right) \text{ for all } t>0 \text{ , } \tau \in(0,1).
 \end{equation}
  Let $u\in W_0^KL_A(\Omega)$ with $[u]_{K,A}<1$.
  Using the definition of Luxemburg-norm and the relation $(\ref{36})$, we deduce 
     $$
     \begin{aligned}
     \displaystyle\int_{Q}A\left( (u(x)-u(y))K(x,y)\right) dxdy&=\displaystyle\int_{Q}A\left([u]_{K,A} \dfrac{(u(x)-u(y))K(x,y)}{[u]_{K,A}} \right) dxdy\\
     &\leqslant [u]^{p_0}_{K,A}\displaystyle\int_{Q}A\left( \dfrac{(u(x)-u(y))K(x,y)}{[u]_{K,A}} \right) dxdy\\
     &\leqslant [u]^{p_0}_{K,A}.
     \end{aligned}
     $$
     
     Finally, we show that $\phi(u) \geqslant [u]^{p_0}_{K,A} \text{   for all  }  u\in W^K_0L_A(\Omega) \text{ with }[u]_{K,A}<1$. Similar techniques as those used in the proof of relation $(\ref{34})$ and $(\ref{35})$, we have 
     \begin{equation}\label{37}
     A(t)\geqslant \tau^{p^0}A\left( \dfrac{t}{\tau} \right) \text{ for all } t>0 \text{ , } \tau \in(0,1).
     \end{equation}
      Let $u\in W_0^KL_A(\Omega)$ with $[u]_{K,A}<1$ and $\beta\in (1,[u]_{K,A})$, so by $ (\ref{37})$ we find 
 \begin{equation}\label{38}
 \displaystyle\int_{Q}A\left( (u(x)-u(y))K(x,y)\right) dxdy\geqslant \beta^{p^0}\int_{Q}A\left( \dfrac{(u(x)-u(y))K(x,y)}{\beta} \right) dxdy.
 \end{equation}
      We define $v(x)=\dfrac{u(x)}{\beta}$ for all $x\in \Omega$, we have $[v]_{K,A}=\dfrac{[u]_{K,A}}{\beta}>1$. Using the relation $(\ref{32})$ we find
      \begin{equation}\label{39}
      \displaystyle\int_{Q}A\left( (v(x)-v(y))K(x,y)\right) dxdy>[v]^{p_0}_{K,A}>1,
      \end{equation}
      by $(\ref{38})$ and $(\ref{39})$ we obtain 
      $$\displaystyle\int_{Q}A\left( (u(x)-u(y))K(x,y)\right) dxdy\geqslant \beta^{p^0}.$$
      Letting $\beta\nearrow [u]_{K,A}$, we deduce that relation $(\ref{33})$ hold true.\\
     \hspace*{15cm$\Box$ }        
\begin{thm}
Let $\Omega$ be an open subset of $\R^N$, and let $K$ as defined in $(K_0)$. The space $W^K_0L_A(\Omega)$ is a Banach space and a separable (resp. reflexive) space if and only if $(A,\Omega)$ is $\Delta$-regular (resp. $(A,\Omega)$ and $(\overline{A},\Omega)$ are $\Delta$-regular). Furthermore
if   $(A,\Omega)$ is $\Delta$-regular and $A(\sqrt{t})$ is convex, then  the space $W^K_0L_A(\Omega)$ is uniformly convex.
\end{thm} 
\noindent \textbf{Proof}.
Let $\left\lbrace u_n\right\rbrace $ be a Cauchy sequence in $W^K_0L_A(\Omega)$. Thus, for any $\varepsilon>0,$ there exist $\mu_\varepsilon>$ such that $n,m\geqslant \mu_\varepsilon$ :
\begin{equation}\label{17}
||u_n-u||_A\leqslant ||u_n-u||_{K,A}<\varepsilon,
\end{equation}
by the completeness of $L_A(\Omega)$, there exist $u\in L_A(\Omega) $ such that $u_n\longrightarrow u$ strongly in $L_A(\Omega)$. Since $u_n=0$ in $\R^N\setminus\Omega$, we define $u=0$ in $\R^N\setminus \Omega$, then $u_n\longrightarrow u$ strongly in $L_A(\R^N)$. So there exist a subsequence $\left\lbrace u_{n_k}\right\rbrace$ in $W^K_0L_A(\Omega)$ such that $u_{n_k}\longrightarrow u$ a.e. in $\R^N$, then by Fatou lemma we get 
$$\int_{Q}  A\left( (u(x)- u(y))K(x,y)\right) dxdy\leqslant \liminf\int_{Q}  A\left( (u_{n_k}(x)- u_{n_k}(y))K(x,y)\right) dxdy<\infty. $$
Thus $u\in W^K_0L_A(\Omega)$. Let $n\geqslant \mu_\varepsilon$, by the second inequality in $(\ref{17})$ and the Fatou lemma we get :
$$||u_n-u||_{K,A}\leqslant \liminf ||u_n-u_{n_k}||_{K,A}\leqslant\varepsilon,$$
that is $u_n\longrightarrow u$ strongly in $W^K_0L_A(\Omega)$.\\
Next, we prove that  $W^K_0L_A(\Omega)$ is a separable, reflexive space and uniformly convex.

We consider the operator 
$$T : W^K_0L_A(\Omega) \longrightarrow  L_A(Q,dxdy).$$ 
 Clearly, T is an isometry. Since $ L_A(Q,dxdy)$ is a reflexive separable and uniformly convex space (see \cite{1,10.}), then $W_0^KL_A(\Omega)$ is also a reflexive separable  and uniformly convex space.\\
 \hspace*{15cm $\Box$ }
      \begin{thm}\label{3.3..}
  Let $A$ be an $N$-function, and $K$ as defined in $(K_0)$. Let $\Omega$  be a bounded open
   subset of  $\R^N$ with $C^{0,1}$-regularity 
     and bounded boundary. If $(\ref{5})$ and $(\ref{6.})$  hold, then 
  \begin{equation}
   W^K{L_A}(\Omega)\hookrightarrow L_{A_*}(\Omega),
  \end{equation}
  and the embedding 
   \begin{equation}
      W^K{L_A}(\Omega)\hookrightarrow L_{B}(\Omega),
     \end{equation}
              is compact for all $B\prec\prec A_*$.
 \end{thm}
 \noindent \textbf{Proof}. Let $u\in W^KL_A(\Omega)$, so by Lemma $\ref{3.2}$, we get $u\in W^sL_A(\Omega)$ and 
 $$ ||u||_{s,A}\leqslant C_1||u||_{K,A}.$$
 On the other hand by Theorem $\ref{2.1}$ we have 
 $$||u||_{A_*}\leqslant C_2||u||_{s,A},$$
 this implies that 
  $$||u||_{A_*}\leqslant C||u||_{K,A}$$
  where $C=C(s,N,k_0,\Omega)$.\\
  By lemma $\ref{3.2}$ and Theorem $\ref{2.2.}$, we have 
      $$  W^K{L_A}(\Omega) \hookrightarrow W^sL_A(\Omega) \hookrightarrow L_{B}(\Omega).$$
       The latter embedding being compact, so we have the desired result.\\
        \hspace*{15cm $\Box$ }
 
 Now, let $A$ be an N-function with $(A,\Omega)$ is $\Delta$-regular. Since $\lim\limits_{t\rightarrow 0}\frac{A(t)}{t}=0$, so there exists $\alpha>0$ such that $A(t)\leqslant t$ for all $t\leqslant \alpha$. We set $\beta=\min\left\lbrace 1,\alpha\right\rbrace $, and we introduce the function $A_1$ as, 
      
       \begin{equation}\label{22.}
       A_1(t)= \left \{
        \begin{array}{clclc}
       \frac{A(\beta)}{\beta} t\hspace{1cm}  & if & t\leqslant \beta,\\\\
          A(t) \hspace*{1cm} & if & t>\beta.
        \end{array}
        \right .  
        \end{equation}
    $A_1$ is a convex, continuous, nondecreasing, finite valued  function, with  $A_1(0)=0$ and  {\small$\lim\limits_{t\rightarrow +\infty} A_1(t)=+\infty$}. $A_1$ is called a Young function (see. \cite{33.} ).\\ 
     For a given domain $\Omega$ in $\R^N$, we define the space $L_{A_1}(\Omega)$ as follows,
     $$L_{A_1}(\Omega) =\Bigg\{u : \Omega\rightarrow \R :    \int_{\Omega}A_1( \lambda u(x))dx < \infty \text{ for some } \lambda>0\Bigg\},$$ 
    this space is equipped with the norm,    
    
       \begin{equation}
        ||u||_{A_1}=\inf \Bigg\{\lambda>0  : \int_{\Omega}A_1\left(  \dfrac{u(x)}{\lambda}\right) dx \leqslant 1\Bigg\}.
       \end{equation}
       \begin{lem}\label{3.3.} (cf. \cite{3})
       Let $\Omega$ be a bounded open subset of $\R^N$ and let $s\in (0,1)$. Let $A$ be an N-function and $A_1$ as defined in $(\ref{22.})$. Then,
       \begin{enumerate}
       \item $L_{A_1}(\Omega)=L_{A}(\Omega)$.
       \item The norm $||.||_A$ and $||.||_{A_1}$ are equivalent.
       \end{enumerate}
       \end{lem}  
 \begin{rem}
 Let $\Omega$ be an open subset of $\R^N$ and let $s\in (0,1)$. Let $A_1$ as defined in $(\ref{22.})$. Then we define the space $W^sL_{A_1}(\Omega)$ by,
 \begingroup\makeatletter\def\f@size{9.5}\check@mathfonts
 \begin{equation}
 W^s{L_{A_1}}(\Omega)=\Bigg\{u\in L_{A_1}(\Omega) : \int_{\Omega} \int_{\Omega} A_1\left( \dfrac{\lambda( u(x)- u(y))}{|x-y|^sA_1^{-1}(|x-y|^N)}\right) dxdy< \infty \text{ for some } \lambda>0 \Bigg\},
 \end{equation}\endgroup
 
 which equipped with the norm
 $$||u||_{s,A_1}=||u||_{A_1}+[u]_{s,A_1}$$
 where, $$[u]_{s,A_1}=\inf \Bigg\{\lambda >0 :  \int_{\Omega} \int_{\Omega} A_1\left( \dfrac{u(x)- u(y)}{\lambda|x-y|^sA_1^{-1}(|x-y|^N)}\right) dxdy\leqslant 1 \Bigg\}.$$
 If  $\Omega$ is a bounded open subset of $\R^N$, then by lemma $\ref{3.3.}$, we have $W^sL_{A_1}(\Omega)=W^sL_{A}(\Omega)$ and the norm $||.||_{s,A}$ and $||.||_{s,A_1}$ are equivalent.
 \end{rem}
 \begin{lem}\label{3.4.}
 Let $\Omega$ be a bounded open subset of $\R^N$ and let $s\in (0,1)$. Let $A_1$ as defined in $(\ref{22.})$. Then the space $W^sL_{A_1}(\Omega)$ continuously embedded in $W^{s,p_0}(\Omega)$ where $p_0>1$ is as defined in $(\ref{12})$. Therefore $W^sL_{A}(\Omega)$ continuously embedded in $W^{s,p_0}(\Omega)$.
 \end{lem}
  \noindent \textbf{Proof}. By definition of  function $A_1$ we get 
\begin{equation}\label{23}
  A(t)\leqslant \gamma A_1(t) \text{ for all } t>0,
\end{equation}
where $\gamma=\max\left\lbrace 1,\frac{\beta}{A(\beta)}\right\rbrace $, and since $A$ satisfies the $\Delta_2$-condition we have 
$$p_0:=\inf_{s\geqslant 0}\dfrac{sa(s)}{A(s)}>1$$
this fact implies that 
$$\dfrac{p_0}{t}\leqslant \dfrac{a(t)}{A(t)} \hspace*{1cm} \text{ for all } t>0,$$
so
$$p_0 (\log(t))'\leqslant (\log(A(t)))' \hspace*{1cm} \text{ for all } t>0$$
therefore
\begin{equation}\label{24}
|t|^{p_0}\leqslant C_1A(t) \hspace*{1cm} \text{ for all } t>0
\end{equation}
we combining $(\ref{23})$ and $(\ref{24})$ we obtain 
\begin{equation}\label{25}
|t|^{p_0}\leqslant C_1A_1(t) \hspace*{1cm} \text{ for all } t>0.
\end{equation}
Let $u\in W^sL_A(\Omega)$, by $(\ref{25})$ we have 
$$\int_{\Omega} |u|^{p_0}dx\leqslant C_1\int_{\Omega}A_1(u)dx $$
then 
\begin{equation}\label{26}
||u||_{L^{p_0}}\leqslant C_1 ||u||_{A_1}.
\end{equation}
where $C_1>0$ and it is possibly different step by step.
On the other hand 
 {\small$$
     \begin{aligned}
  \int_{\Omega}\int_{\Omega} \dfrac{|u(x)-u(y)|^{p_0}}{|x-y|^{p_0s+N}}dxdy&=  \int_{\Omega}\int_{\Omega\cap\left\lbrace  |x-y|^N\leqslant \beta\right\rbrace } \dfrac{|u(x)-u(y)|^{p_0}}{|x-y|^{p_0s+N}}dxdy\\
  & \hspace*{1cm}+  \int_{\Omega}\int_{\Omega\cap\left\lbrace |x-y|^N> \beta\right\rbrace } \dfrac{|u(x)-u(y)|^{p_0}}{|x-y|^{p_0s+N}}dxdy\\
  &= I_1+I_2.
     \end{aligned}
$$}
By definition of $A_1$ and estimation $(\ref{25})$, we have
\begingroup\makeatletter\def\f@size{8}\check@mathfonts
\begin{equation}\label{27.}
     \begin{aligned}
  I_1=\int_{\Omega}\int_{\Omega\cap\left\lbrace  |x-y|^N\leqslant \beta\right\rbrace } \dfrac{|u(x)-u(y)|^{p_0}}{|x-y|^{p_0s+N}} & dxdy=\frac{A(\beta)}{\beta}\int_{\Omega}\int_{\Omega\cap\left\lbrace  |x-y|^N\leqslant \beta\right\rbrace } \dfrac{|u(x)-u(y)|^{p_0}}{|x-y|^{sp_0}A_1^{-1}(|x-y|^N)}dxdy\\
 & \leqslant \left( \frac{A(\beta)}{\beta}\right) ^{p_0}\int_{\Omega}\int_{\Omega\cap\left\lbrace  |x-y|^N\leqslant \beta\right\rbrace } \dfrac{|u(x)-u(y)|^{p_0}}{|x-y|^{sp_0}(A_1^{-1}(|x-y|^N))^{p_0}}dxdy\\
 &\leqslant C_1\left( \frac{A(\beta)}{\beta}\right) ^{p_0}\int_{\Omega} \int_{\Omega} A_1\left( \dfrac{u(x)- u(y)}{\lambda|x-y|^sA_1^{-1}(|x-y|^N)}\right) dxdy
     \end{aligned}
\end{equation}
\endgroup
and 
\begingroup\makeatletter\def\f@size{10}\check@mathfonts
$$
     \begin{aligned}
  I_2&=\int_{\Omega}\int_{\Omega\cap\left\lbrace  |x-y|^N> \beta\right\rbrace } \dfrac{|u(x)-u(y)|^{p_0}}{|x-y|^{sp_0+N}}dxdy\\
  &=\int_{\Omega}\int_{\Omega\cap\left\lbrace  |x-y|^N> \beta\right\rbrace } \dfrac{|u(x)-u(y)|^{p_0}}{|x-y|^{sp_0}(A_1^{-1}(|x-y|^N))^{p_0}}\dfrac{(A_1^{-1}(|x-y|^N))^{p_0}}{|x-y|^N}dxdy\\
 & \leqslant \sup_{\Omega \times \Omega \cap\left\lbrace  |x-y|^N> \beta\right\rbrace}\dfrac{(A_1^{-1}(|x-y|^N))^{p_0}}{|x-y|^N}\int_{\Omega}\int_{\Omega\cap\left\lbrace  |x-y|> \beta\right\rbrace } \dfrac{|u(x)-u(y)|^{p_0}}{|x-y|^{sp_0}(A_1^{-1}(|x-y|^N))^{p_0}}dxdy,\\
 &\leqslant \sup_{\Omega \times \Omega \cap\left\lbrace  |x-y|^N> \beta\right\rbrace}\dfrac{(A_1^{-1}(|x-y|^N))^{p_0}}{|x-y|^N}C_1\int_{\Omega} \int_{\Omega} A_1\left( \dfrac{u(x)- u(y)}{\lambda|x-y|^sA_1^{-1}(|x-y|^N)}\right) dxdy.\\
     \end{aligned}
$$\endgroup
Since $A^{-1}_1(t)$ is continuous for all $t>\beta$ and $\Omega$ is bounded so,
$$\sup_{\Omega \times \Omega \cap\left\lbrace  |x-y|^N> \beta\right\rbrace}\dfrac{A_1^{-1}(|x-y|^N)}{|x-y|^N}=C_2< \infty,$$
therefore
\begin{equation}\label{28}
I_2\leqslant C_1C_2\int_{\Omega} \int_{\Omega} A_1\left( \dfrac{u(x)- u(y)}{\lambda|x-y|^sA_1^{-1}(|x-y|^N)}\right) dxdy,
\end{equation}
by combining $(\ref{27.})$ and $(\ref{28})$ we obtain that 
$$\int_{\Omega}\int_{\Omega} \dfrac{|u(x)-u(y)|^{p_0}}{|x-y|^{p_0s+N}}dxdy\leqslant C \int_{\Omega} \int_{\Omega} A_1\left( \dfrac{u(x)- u(y)}{\lambda|x-y|^sA_1^{-1}(|x-y|^N)}\right) dxdy,$$
where $C_3=C_1\left( \dfrac{A(\beta)}{\beta}\right) ^{p_0}+C_1C_2$. Then 
\begin{equation}\label{29}
[u]_{s,p_0}\leqslant C_3 [u]_{s,A_1}.
\end{equation}
Combining $(\ref{26})$ and $(\ref{29})$, we have
$$ ||u||_{s,p_0}\leqslant C ||u||_{s,A_1},$$
Where $C=C_1+C_3$.\\
        \hspace*{15cm $\Box$ }
\begin{coro}\label{3.1.}
Let $\Omega$ be a bounded open subset of $\R^N$, let $s\in (0,1)$ and let $A$ be an $N$-function satisfies the global $\Delta_2$-condition. Let $K$ as defined in $(K_0)$.\\
 If $sp_0<N$,
        then 
       $$W^KL_A(\Omega) \hookrightarrow L^{p^*_0}(\Omega),$$
     where $
       p^*_0= \frac{Np_0}{N-sp_0}.
             $  And the embedding  
    $$W^KL_A(\Omega) \hookrightarrow L^{q}(\Omega),$$
    is compact for all $q\in [1,p_0^*)$.\\ If $sp_0=N$, then embedding    $$W^KL_A(\Omega) \hookrightarrow L^{q}(\Omega),$$ is compact, for all $q\in [1,+\infty)$.\\
    If $sp>N$, then the embedding
    $$W^KL_A(\Omega) \hookrightarrow L^{\infty}(\Omega),$$ is compact.
\end{coro}
\section{Existence results}
In this section,  we prove the existence of a weak solutions for a fractional Kirchhoff type problem in
fractional Orlicz-Sobolev spaces, by means of  the direct method in calculus variations.
For this, we suppose that the Kirchhoff function
$M: [0,\infty) \longrightarrow (0,\infty)$ is a continuous function satisfying the following condition:

There exists $m_0>0$ such that :
$$ (M_0) \hspace{4.2cm} M(t)\geqslant m_0 \text{  for all  } t\in [0,\infty). \hspace{7cm} $$

Also, we assume that $f: \Omega\times \R \longrightarrow \R$ is a Carath\'eodory function satisfying the following conditions :\\
there exists $\theta_1>0$ and $1 < q < p_0^*$  such that

$$ (f_1) \hspace{4.2cm}|f(x,t)|\leqslant \theta_1 (1+|t|^{q-1}) \text{  a.e.  } (x,t)\in \Omega\times \R^N, \hspace{7cm} $$
there exist $\theta_2>0$ and an open bounded set $\Omega_0\subset \Omega$ such that

$$ (f_2) \hspace{4.2cm}|f(x,t)|\geqslant \theta_2 |t|^{q-1} \text{  a.e.  } (x,t)\in \Omega_0\times \R^N. \hspace{7cm} $$
\begin{defini}
We say that $u\in W^K_0L_A(\Omega)$ is a weak solution of problem $(D_{K,A})$ if 
$$M\left( \int_{\R^N}A(h_{x,y}(u))dxdy\right) \int_{\R^{2N}}a(h_{x,y}(u)\left( h_{x,y}(v)\right) dxdy=\int_{\Omega}f(x,u)vdx,$$
for all $v\in W^K_0L_A(\Omega)$, where $h_{x,y}(u)=(u(x)-u(y))K(x,y)$.
\end{defini}
\begin{thm}\label{4.1.}
Let $K  : \R^N\times\R^N\longrightarrow (0,\infty)$ be a function satisfy $(K_0)$, and $A$  be an $N$- function satisfies $(a_1)$ and $(a_2)$, suppose that $M$ satisfy $(M_0)$, and $f$ satisfies $(f_1)$ and $(f_2)$ if $1<q<p_0$, then the problem $(D_{K,A})$ has a nontrivial weak solution in $W^K_0L_A(\Omega)$. 
\end{thm}
\begin{coro}
Let $K  : \R^N\times\R^N\longrightarrow (0,\infty)$ be a function satisfies $(K_0)$, and $A$  be an $N$- function satisfies $(a_1)$ and $(a_2)$, suppose that $M$ satisfies $(M_0)$, and $f$ satisfy $(f_1)$ and $(f_2)$, if $q=p_0$ and $\theta_1<(m_0\lambda_1)/2$, where $\lambda_1$ is the first eigenvalue of $\mathcal{L}^K_A u$ define by
$$\lambda_1=\inf\limits_{u\in W^k_0L_A(\Omega)\setminus \left\lbrace 0\right\rbrace }\dfrac{||u||^{p_0}_{K,A}}{||u||_A^{p_0}},$$
 then the problem $(D_{K,A})$ has a nontrivial weak solution in $W^K_0L_A(\Omega)$. 
\end{coro}
For $u\in W^K_0L_A(\Omega)$, we define 
$$J(u)=\widehat{M}\left( \int_{Q}A((u(x)-u(y))K(x,y))dxdy\right),$$
$$ H(u)=\int_{\Omega}F(x,u)dx \text{ and } I(u)=J(u)-H(u),$$
where $\widehat{M}(t)=\displaystyle\int_{0}^{t}M(\tau)d\tau$ and $F(x,t)=\displaystyle\int_{0}^{t}f(x,t)d\tau$.
Obviously the energy functional $I : W^K_0L_A(\Omega)\longrightarrow \R$ associated with problem $(D_{K,A})$ is well defined.
\begin{lem}\label{4.1}
If $f$ satisfies assumption $(f_1)$, then the functional $H\in C^1(W^K_0L_A,\R)$ and 
$$<H'(u),v>=\int_{\Omega}f(x,u)vdx \text{ for all } u,v \in W^K_0L_A(\Omega).$$
\end{lem}
\noindent \textbf{Proof}. By corollary $\ref{3.1.}$, the proof of this lemma is similar to proof of lemma 3.1 in \cite{111}.\\
        \hspace*{15cm $\Box$ }
        \begin{lem}\label{4.2}
Let $(M_0)$ hold, then the function $J\in C^1(W^K_0L_A,\R)$ and
$$<J'(u),v>=M\left( \int_{Q}A(h_{x,y}(u))dxdy\right) \int_{Q}a(h_{x,y}(u)) h_{x,y}(v) dxdy$$
for all $u,v \in W^K_0L_A(\Omega)$. Moreover, for   each $u \in W^K_0L_A(\Omega)$, $J'(u) \in (W^K_0L_A(\Omega))^*$ where $(W^K_0L_A(\Omega))^*$, denotes the dual space of $W^K_0L_A(\Omega)$.

        \end{lem}
        \noindent \textbf{Proof}.
   First, it is easy to see that 
\begin{equation}\label{30}
<J'(u),v>=M(\int_{Q}A(h_{x,y}(u)dxdy)\int_{Q}a(h_{x,y}(u)) h_{x,y}(v)dxdy
\end{equation}
 for all $u,v \in W^K_0L_A(\Omega)$. It follows from $(\ref{30})$  that for each $u \in W^K_0L_A(\Omega)$, $J'(u) \in (W^K_0L_A(\Omega))^*$.
 
 Next, we prove that $J\in C^1(W^K_0L_A(\Omega),\R)$. Let $\left\lbrace u_n\right\rbrace \subset W^K_0L_A(\Omega)$ with $u_n\longrightarrow u$ strongly in $W^K_0L_A(\Omega)$,
  for $v\in W^K_0L_A(\Omega)$  we have $h_{x,y}(v)\in L_{A}(Q,dxdy)$ 
                 and by H\"older's inequality,
                 \begingroup\makeatletter\def\f@size{11}\check@mathfonts       
          $$     
      \begin{aligned}
     \left| \int_{Q} ( a(h_{x,y}(u_n))
     -a(h_{x,y}(u))) h_{x,y}(v)dxdy\right|&
     \leqslant 2\left|\left| a(h_{x,y}(u_n))
       -a(h_{x,y}(u))\right| \right|_{L_{\overline{A}}} \left| \left| h_{x,y}(v)\right| \right|_{L_{A}}\\
       &\leqslant 2 \left|\left| a(h_{x,y}(u_n))
           -a(h_{x,y}(u))\right| \right|_{L_{\overline{A}}} \left| \left|h_{x,y}(v)\right| \right|_{K,A}.
   \end{aligned}  
   $$ \endgroup
   
   On the other hand, $u_n \rightarrow u$ in $W^K_0L_A(\Omega)$, then $v_n:=h_{x,y}(u_n)\longrightarrow v:=h_{x,y}(u)$ in $L_A(Q)$, so by dominated convergence theorem, there exists a subsequence $\left\lbrace v_{n_k}\right\rbrace $ and a function $h$ in $L_A(Q)$ such that 
    $$|a(h_{x,y}(u_{n_k}))|\leqslant |a(h)| \in L_{\overline{A}}(Q) \text{ a.e in  } Q,$$
    and
    $$a(h_{x,y}(u_{n_k}))\longrightarrow a(h_{x,y}(u))   \text{ a.e in } Q.$$ 
    Then by dominated convergence theorem we obtain that
    $$\sup_{||v||_{K,A}\leqslant 1} \Bigg| \int_{Q} ( a(h_{x,y}(u_n))
         -a(h_{x,y}(u))) h_{x,y}(v)dxdy\Bigg|\longrightarrow 0.$$
         Moreover, by the continuity of $M$, we have 
         $$M\left( \int_{\R^N}A(h_{x,y}(u_n))dxdy\right) \longrightarrow M\left( \int_{\R^N}A(h_{x,y}(u))dxdy\right). $$
   \hspace*{15cm$\Box$ }
 
 Combining lemma $\ref{4.1}$ and lemma $\ref{4.2}$, we get $I \in C^1(W^K_0L_A(\Omega),\R)$ and 
 {\small$$<I'(u),v>=M\left( \int_{Q}A(h_{x,y}(u))dxdy\right) \int_{Q}a(h_{x,y}(u)\left( h_{x,y}(v)\right) dxdy- \int_{\Omega}f(x,u)vdx$$}
   for all  $u,v \in W^K_0L_A(\Omega)$.\\

   \begin{lem}\label{4.3}
   Let $(M_0)$ and $(f_1)$ be satisfied, then the functional $I\in C^1(W^K_0L_A(\Omega), \R)$ is weakly lower semi-continuous.
   \end{lem}
   \noindent \textbf{Proof}. First, note that the map :
   $$ u\longmapsto \int_{Q}A\left( (u(x)-u(y))K(x,y)\right) dxdy $$ 
   is lower semi-continuous for the weak topology of $W^K_0L_A(\Omega)$. Indeed, we define a functional $\psi :  W^K_0L_A(\Omega) \longrightarrow \R $ as 
   $$\psi(u)=\int_{Q}A\left( (u(x)-u(y))K(x,y)\right) dxdy,$$
   similar to Lemma $\ref{4.2}$, we obtain $\psi \in C^1(W^K_0L_A(\Omega), \R^N)$ and
   $$< \psi'(u),v>=\int_{Q}a\left( (u(x)-u(y))K(x,y)\right)(v(x)-v(y))K(x,y)dxdy $$
      for all  $u,v \in W^K_0L_A(\Omega)$.\\
       On the other hand, since $A$ is convex so $\psi$ is also convex.
      	Now, let $\left\lbrace u_n\right\rbrace \subset W^K_0L_A(\Omega)$ with $u_n\rightharpoonup u$ weakly in $ W^K_0L_A(\Omega)$, then by convexity of $\psi$ we have 
      	$$\psi(u_n)-\psi(u)\geqslant <\psi'(u),u_n-u>,$$ 
      	hence, we obtain $\psi(u)\leqslant \liminf\psi(u_n)$, that is, the map
 $$ u\longmapsto \int_{Q}A\left( (u(x)-u(y))K(x,y)\right) dxdy $$      	
   is lower semi-continuous.
   
  Let $u_n\rightharpoonup u$ weakly in $W^k_0L_A(\Omega)$, so by corollary $\ref{3.1.}$, $u_n\longrightarrow u$ in $L^q(\Omega)$ for all $q\in (p_0,p_0^*)$. Without loss of generality, we assume that $u_n\longrightarrow u$ a.e. in $\Omega$. Assumption $(f_1)$ implies that 
  $$F(x,t)\leqslant 2\theta_1 (|t|^q+1).$$
  Thus, for any measurable subset $U\subset \Omega$,
  $$\int_{U}|F(x,u_n)|dx\leqslant 2\theta_1 \int_{U}|u_n|^qdx+2\theta_1|U|.$$
  By H\"older inequality and corollary $\ref{3.1.}$ , we have,

  \begin{equation}\label{37.}
         \begin{aligned}
        \int_{U}|F(x,u_n)|dx&\leqslant 2\theta_1 ||u_n^q||_{L^{\frac{p^*_0}{q}}}||1||_{L^{\frac{p^*_0}{p_0^*-q}}}+2\theta_1|U|\\
        &\leqslant 2\theta_1C ||u_n||^q_{K,A}|U|^{\frac{p_0^*-q}{p_0^*}}+2\theta_1|U|.
      \end{aligned}  
  \end{equation} 
    
  It follows from $(\ref{37.})$ that the sequence $\left\lbrace |F(x,u_n)-F(x,u)|\right\rbrace $ is uniformly bounded and equi-integrable
  in $L^1(\Omega)$. The Vitali Convergence Theorem (see \cite{013.}) implies
  $$\lim\limits_{n\rightarrow n}\int_{\Omega}|F(x,u_n)-F(x,u)|dx=0,$$
   so
   $$\lim\limits_{n\rightarrow \infty}\int_{\Omega}F(x,u_n)dx=\int_{\Omega}F(x,u)dx.$$
   Thus, the functional $H$ is weakly continuous. Further, we get that $I$ is weakly lower semi-continuous.\\
     \hspace*{15cm$\Box$ } 
        \noindent \textbf{Proof of theorem $\ref{4.1.}$}.  
   From assumptions $(M_0)$, $(f_1)$  and proposition  $\ref{pro3}$, we have 
   $$
   \begin{aligned}
          I(u)&=\widehat{M}\left( \int_{Q} A\left( (u(x)-u(y))K(x,y)\right) dxdy\right) -\int_{\Omega}F(x,u)dx\\
          &\geqslant m_0\int_{Q} A\left( (u(x)-u(y))K(x,y)\right) dxdy-2\theta_1 \int_{\Omega} |u|^qdx-2\theta_1|\Omega|\\
          &\geqslant m_0[u]_{K,A}^{p_0} -2\theta_1 C [u]_{K,A}^q-2\theta_1|\Omega|,
         \end{aligned}
         $$
since $p_0>q$, so we have $I(u)\longrightarrow \infty$ as $[u]_{K,A}\longrightarrow \infty$, by Lemma $\ref{4.3}$ $I$ is weakly lower semi-continuous on $W^K_0L_A(\Omega)$, then by theorem \ref{2.2} functional I has a minimum point $u_0$ in $W^K_0L_A(\Omega)$ and $u_0$ is a weakly solution of problem $(D_{K,A})$.

Next we need to verify that $u_0$ is nontrivial. Let $x_0\in \Omega_0$, $0 < R < 1$ satisfy $B_{2R}(x_0)\subset \Omega_0$, where $B_{2R}(x_0)$
is the ball of radius $2R$ with center at the point $x_0$ in $\R^N$. Let $\varphi\in C_0^\infty(B_{2R}(x_0))$ satisfies $0\leqslant \varphi \leqslant 1$ and
$\varphi \equiv 1$  in $B_{2R}(x_0)$. Lemma $\ref{3.1..}$ implies that $$||u||_{K,A}<\infty.$$ Then for $0 < t < 1$, by the mean value theorem
and $(f_2)$, we have

 $$
   \begin{aligned}
          I(t\varphi)&=\widehat{M}\left( \int_{Q} A\left( (t\varphi(x)-t\varphi(y))K(x,y)\right) dxdy\right) -\int_{\Omega}F(x,t\varphi)dx\\
          &\leqslant \widehat{M}\left(||t\varphi||^{p^0}_{K,A}\right) -\int_{\Omega}F(x,t\varphi)dx \\
          &\leqslant  \int_{0}^{||t\varphi||^{p^0}_{K,A}} M(\tau)d\tau-\int_{\Omega_0}\dfrac{\theta_2}{q} |t\varphi|^{q}dx\\
          &\leqslant M(v)||\varphi||^{p^0}_{K,A}t^{p^0}-\dfrac{\theta_2}{q}t^q\int_{\Omega_0} |\varphi|^{q}dx\\
          &\leqslant Ct^{p^0}-\dfrac{\theta_2}{q}t^q\int_{\Omega_0} |\varphi|^{q}dx,
         \end{aligned}
         $$
where $v\in [0,||\varphi||^{p^0}_{k,A})$ and $C$ is a positive constant. Since $p^0 > q$ and $\displaystyle\int_{\Omega_0} |\varphi|^{q}dx>0,$ we have $I(t_0\psi)<0$ for
$t_0\in(0,t)$ sufficiently small. Hence, the critical point $u_0$ of functional $I$ satisfies $I(u_0)\leqslant I(t_0\psi)< 0=I(0)$,
that is $u_0\neq 0$. \\
        \hspace*{15cm $\Box$ }\\
 \noindent \textbf{Proof  of corollary }. In view of the proof of theorem $\ref{4.1.}$, we only need to check that $I(u)\longrightarrow \infty$ as $||u||_{K,A}\rightarrow \infty$. Since
 $p = q$ and $\theta_1< (m_0\lambda_1)/(2)$, by assumption $(f_1)$ and the definition of first eigenvalue of $\mathcal{L}^K_A u$ , we have  
 
  $$
    \begin{aligned}
           I(u)&=\widehat{M}\left( \int_{Q} A\left( (u(x)-u(y))K(x,y)\right) dxdy\right) -\int_{\Omega}F(x,u)dx\\
           &\geqslant m_0\int_{Q} A\left( (u(x)-u(y))K(x,y)\right) dxdy-2\theta_1 \int_{\Omega} |u|^{p_0}dx-2\theta_1|\Omega|\\
           &\geqslant m_0[u]_{K,A}^{p_0} -2\theta_1 \frac{1}{\lambda_1} [u]_{K,A}^{p_0}-2\theta_1|\Omega|,\\
           & =\left( m_0 -2\theta_1\frac{1}{\lambda_1}\right) [u]_{K,A}^{p_0}-2\theta_1|\Omega|.
          \end{aligned}
          $$
          So we have $I(u)\longrightarrow \infty$ as $[u]_{K,A}\longrightarrow \infty$.\\
               \hspace*{15cm$\Box$ }

      \end{document}